\def\draft{n}
\documentclass{amsart}
\usepackage{fullpage,amssymb,epic,eepic,epsfig,amsbsy,amsmath}

\theoremstyle{plain}

\newtheorem{theorem}{Theorem}
\newtheorem{proposition}{Proposition}[section]
\newtheorem{lemma}[proposition]{Lemma}
\newtheorem{corollary}[proposition]{Corollary}

\theoremstyle{definition}
\newtheorem{definition}[proposition]{Definition}

\newtheorem{question}{Question}

\theoremstyle{remark}
\newtheorem{example}[proposition]{Example}

\newtheorem{remark}[proposition]{Remark}

\def\printname#1{
        \if\draft y
                \smash{\makebox[0pt]{\hspace{-0.5in}
                        \raisebox{8pt}{\tt\tiny #1}}}
        \fi
}

\newcommand{\psdraw}[2]
         {\begin{array}{c} \hspace{-1.3mm}
        \raisebox{-4pt}{\epsfig{figure=draws/#1.eps,width=#2}}
        \hspace{-1.9mm}\end{array}}

\newlength{\standardunitlength}
\catcode`\@=11
\long\def\@makecaption#1#2{%
     \vskip 10pt

\setbox\@tempboxa\hbox{
       \small\sf{\bfcaptionfont #1. }\ignorespaces #2}%
     \ifdim \wd\@tempboxa >\captionwidth {%
         \rightskip=\@captionmargin\leftskip=\@captionmargin
         \unhbox\@tempboxa\par}%
       \else
         \hbox to\hsize{\hfil\box\@tempboxa\hfil}%
     \fi}
\font\bfcaptionfont=cmssbx10 scaled \magstephalf
\newdimen\@captionmargin\@captionmargin=2\parindent
\newdimen\captionwidth\captionwidth=\hsize
\catcode`\@=12

\def\lbl#1{\label{#1}\printname{#1}}

\def\BZ{\mathbb Z}

\def\BQ{\mathbb Q}
\def\BR{\mathbb R}
\def\BC{\mathbb C}

\def\BT{\mathbb T}

\def\A{\mathcal A}
\def\B{\mathcal B}
\def\C{\mathcal C}
\def\D{\Delta}

\def\R{\mathcal R}

\def\calS{\mathcal S}
\def\calI{\mathcal I}

\def\calP{\mathcal P}
\def\calO{\mathcal O}

\def\a{\alpha}

\def\l{\lambda}

\def\la{\langle}
\def\ra{\rangle}

\def\ti{\widetilde}

\def\longto{\longrightarrow}

\def\pt{\partial}
\def\fg{\mathfrak{g}}
\def\Zrat{Z^{\mathrm{rat}}}
\def\Jrat{J^{\mathrm{rat}}}
\def\Ev{\mathrm{Ev}}
\def\calD{\mathcal{D}}
\def\calX{\mathcal{X}}
\def\Lhat{\hat{\Lambda}}
\def\invlim{\mathrm{invlim}}

\begin{document}


\title[Difference and differential equations for the colored Jones 
function]{Difference and differential equations for the colored Jones 
function}

\author{Stavros Garoufalidis}
\address{School of Mathematics \\
          Georgia Institute of Technology \\
          Atlanta, GA 30332-0160, USA. \\ \newline
         {\tt http://www.math.gatech.edu/$\sim$stavros } }
\email{stavros@math.gatech.edu}

\thanks{Supported in part by the National Science Foundation. \\
\newline
1991 {\em Mathematics Classification.} Primary 57N10. Secondary 57M25.
\newline
{\em Key words and phrases: holonomic function, colored Jones function, 
recursion ideal, peripheral ideal, orthogonal ideal, 
Kauffman bracket skein module, loop expansion, hierarchy of ODE, WKB.} 
}

\date{
July 25, 2005} 


\begin{abstract}
The colored Jones function of a knot is a sequence of Laurent polynomials.
It was shown by TTQ. Le and the author that such sequences are $q$-holonomic,
that is, they satisfy linear $q$-difference equations with coefficients
Laurent polynomials in $q$ and $q^n$. We show from first principles that
$q$-holonomic sequences give rise to modules over a $q$-Weyl ring.
Frohman-Gelca-LoFaro have identified 
the latter ring with the ring of even functions of the quantum
torus, and with the Kauffman bracket skein module of the torus.
Via this identification, we study relations among the orthogonal,
peripheral and recursion ideal of the colored Jones function, introduced
by the above mentioned authors.
In the second part of the paper, we convert the linear $q$-difference
equations of the colored Jones function in terms of a hierarchy of linear
ordinary differential equations for its loop expansion. This conversion
is a version of the WKB method, and may shed some information on the problem
of asymptotics of the colored Jones function of a knot.
\end{abstract}

\maketitle



\section{Introduction}
\lbl{sec.intro}

\subsection{The colored Jones function and its loop expansion}
\lbl{sub.goal}
The {\em colored Jones function} of a knot $K$ in 3-space is a sequence
$$
J_K: \BZ \longto \BZ[q^{\pm/2}]
$$
of Laurent polynomials 
that encodes the {\em Jones polynomial} of a knot and its
parallels; \cite{J,Tu}. Technically, $J_{K,n}$ is the quantum group invariant
using the $n$-dimensional representation of $\mathfrak{sl}_2$ for $n \geq 0$,
normalized by 
$J_{\text{unknot},n}(q)=[n]$ (where 
$[n]=(q^{n/2}-q^{-n/2})/(q^{1/2}-q^{-1/2})$), 
and extended to integer indices by 
$J_{K,n}=-J_{K,-n}$. 

In the spring of 2005, TTQ. Le and the author proved
that the colored Jones function is $q$-{\em holonomic}, \cite{GL1}. 
In in other words, it satisfies a {\em linear $q$-difference equation} 
with coefficients Laurent polynomials in $q$ and $q^n$. 

In \cite{Ro}, Rozansky introduced a {\em loop expansion} of the colored Jones
function. Namely, he associated to a knot $K$ an invariant
$$
\Jrat_K(q,u)= \sum_{k=0}^\infty Q_{K,k}(u)(q-1)^k \in \BQ'(u)[[q-1]]
$$
where $\BQ'(u)[[q-1]]$ is the ring of power series in $q-1$ 
with coefficients in the ring $\BQ'(u)$ of rational functions in $u$ 
which do not have a pole at $u=1$, and 
$$
Q_{K,k}(u)=\frac{P_{K,k}(u)}{\D_K(u)^{2k+1}}
$$
where $P_{K,k}(u) \in \BQ[u^{\pm}]$ and $\D_K(t)$ is the 
{\em Alexander polynomial}
normalized by $\D_K(t)=\D_K(t^{-1})$, $\D_K(1)=1$ and 
$\D_{\text{unknot}}(t)=1$.

The relation between $\Jrat_K$ and $J_K$ is the following equality,
valid in the power series ring $\BQ[[q-1]]$ 
\begin{eqnarray}
\lbl{eq.Jrat}
[n] \Jrat_K(q,q^n)=J_{K,n}(q) \in \BQ[[q-1]],
\end{eqnarray}
for all $n > 0$. Notice that $\D_K(1)=1$, thus $1/\D_K(u) \in \BQ'(u)$
and $1/\D_K(q^n)$ can always
be expanded in power series of $q-1$. Notice moreover that $\Jrat_K$
determines $J_K$ and vice-versa, via the above equation. 

In this paper, we convert linear $q$-difference equations for the 
colored Jones function $J_K$ into a 
{\em hierarchy of linear differential equations} for 
the loop expansion $\Jrat_K$.

Moreover, we study holonomicity of the colored Jones function from the
point of view of quantum field theory, and compare it   
with the skein theory approach of the
colored Jones function initiated by Frohman and Gelca.

The paper was written in the spring of 2003, following stimulating 
conversations with A. Sikora, who kindly explained to the author the work
of Gelca-Frohman-LoFaro and others on the skein theory approach to the
colored Jones function. The author wishes to thank A. Sikora for 
enlightening conversations. The paper remained as a preprint for over two
years. The current version is substantially revised to take into account
the recent developments of the last two years.

\subsection{Holonomic functions and the $q$-Weyl ring $\C$}
\lbl{sub.holonomic}

A holonomic function $f(x)$ in one continuous variable $x$ is one 
that satisfies a nontrivial linear differential equation with polynomial 
coefficients.

In this section, we will see that the notion of holonomicity for a sequence
of Laurent polynomials {\em naturally leads} 
to a $q$-Weyl ring $\C$ defined below.

Holonomicity was introduced by I.N. Bernstein \cite{B1,B2} in relation
to algebraic geometry, $D$-modules and differential Galois theory.
In a stroke of brilliance, Zeilberger noticed that holonomicity can be
applied to verify, in a systematic way, combinatorial 
identities among special functions, \cite{Z}. This was later implemented
on a computer, \cite{WZ, PWZ}.

A key  idea is to study the recursion relations that a function
satisfies, rather than the function itself. This idea leads in a natural way
to noncommutative algebras of operators that act on a function, together with
left ideals of annihilating operators.

To explain this idea concretely, consider the operators $x$ and $\pt$ which 
act on a 
smooth function $f$ defined on $\BR$ (or a distribution, or whatever else can 
be differentiated) by
$$
(xf)(x) = x f(x) \hspace{2cm}
(\pt f)(x) = \frac{\pt}{\pt x} f(x).
$$
Leibnitz's rule $ \pt ( x f)= x \pt(f) +  f$ written in operator form
states that $\pt x = x \pt +1$. The operators $x$ and $\pt$ generate
the {\em Weyl algebra} which is a free noncommutative algebra on $x$ and $\pt$
modulo the two sided ideal $\pt x - x \pt -1$:
$$
\A = \frac{\BC\la x, \pt  \ra}{
(\pt x - x \pt -1)}.
$$ 
The Weyl algebra is nothing but the {\em algebra of differential operators
in one variable with polynomial coefficients}.
Given a function $f$ of one variable, let us define the 
{\em recursion ideal} $\calI_f$ by
$$
\calI_f=\{ P \in \A | \, Pf=0 \}
$$
It is easy to see that $\calI_f$ is a left ideal of $\A$.
Following Zeilberger and Bernstein, we say that $f$ is {\em holonomic} 
iff $\calI_f \neq 0$. In other words, a holonomic function is one that
satisfy a linear differential equation with polynomial coefficients.

A key property of the Weyl algebra $\A$ (shared by its cousins, $\B$ and
$\C$ defined below) is that it is {\em Noetherian}, which implies that
every left ideal is {\em finitely generated}. In particular, a holonomic
function is uniquely determined by a finitely list, namely the generators
of its recursion ideal and a finite set of initial conditions.

The set of holonomic functions is closed under summation and product.
Moreover, holonomicity can be extended to functions of several variables. 
For an excellent exposition of these results, see \cite{Bj}.

Zeilberger expanded the definition of holonomic functions of a continuous
variable to 
{\em discrete functions} $f$ (that is, functions with domain $\BZ$; otherwise
known as bi-infinite sequences) by 
replacing differential
operators by {\em shift} operators. More precisely, consider the operators
$N$ and $E$ which act on a discrete function 
$(f_n)$ by
$$
(N f)_n = n f_n, \quad
(E f)_n = f_{n+1}.
$$
It is easy to see that $EN=NE+E$.
The {\em discrete Weyl algebra} $\B$ is a noncommutative algebra with
presentation
$$
\B= \frac{\BQ\la N^{\pm}, E^{\pm} \ra}{
(EN-NE-E)}.
$$
The field coefficients $\BQ$ are not so important, and neither is
the fact that we allow positive as well as negative powers of $N$ and $E$.
Given a discrete function $f$, one can define the {\em recursion ideal} 
$\calI_f$ in
$\B$ as before. We will call a discrete function $f$ {\em holonomic}
iff the ideal $\calI_f \neq 0$.

In our paper we will consider a $q$-variant of the Weyl algebra. 
Let
\begin{equation}
\lbl{eq.R}
\R=\BZ[q^{\pm/2}].
\end{equation}
Consider the
operators $E$ and $Q$ which act on a discrete function 
$f: \BZ \longto \R$ by:
\begin{equation}
\lbl{eq.EQaction1}
(Q f)_n(q) = q^n f_n(q), \quad
(E f)_n(q) = f_{n+1}(q).
\end{equation}
It is easy to see that $EQ=qQE$.
We define the $q$-{\em Weyl ring} $\C$ to be a noncommutative ring 
with presentation
\begin{equation}
\lbl{eq.qWeyl}
\C= \frac{\R\la Q, E \ra}{
(EQ-qQE)}.
\end{equation}
Given a discrete function $f:\BZ \longto \R$, one can define the left 
ideal $\calI_f$ in $\C$ as before, and call a discrete function $f$ 
$q$-{\em holonomic} iff the ideal $\calI_f \neq 0$.
Concretely, a discrete function $f: \BZ \longto \R$ is $q$-holonomic 
iff there exists a nonzero element $ \sum_{a,b} c_{a,b} E^a Q^b \in \C$ such 
that
\begin{equation}
\lbl{eq.recursion}
\sum_{a,b} c_{a,b} q^{(n+a)b} f_{n+a}(q)=0.
\end{equation}

The sequence
$J: \BZ \longto \R$ 
of Laurent polynomials that we have in mind is the celebrated {\em
colored Jones function}. 

\begin{theorem}
\cite{GL1}
\lbl{thm.GL}
The colored Jones function of every knot in $S^3$ is $q$-holonomic.
\end{theorem}


\subsection{The $q$-torus
and the Kauffman bracket skein module of the torus}
\lbl{sub.2incarnations}

The above section we introduced the $q$-Weyl ring $\C$ to define the
notion of holonomicity of a sequence of Laurent polynomials.

The $q$-Weyl ring is not new to quantum topology. It has already appeared
in the theory of {\em quantum groups}, (see \cite[Chapter IV]{Ka}, \cite{M}), 
under the name: the algebra of functions of the  {\em quantum torus}.

It has also appeared in the skein theory approach of the colored Jones
polynomial, via the Kauffman bracket. Let us review this important
discovery of Frohman, Gelca and Lofaro, \cite{FGL}.
As a side bonus, we can associate two further natural knot invariants:
the quantum peripheral and orthogonal ideals.

\subsection{The quantum peripheral and orthogonal ideals of a knot}
\lbl{sub.results}

The recursion relations \eqref{eq.recursion} for the colored Jones function
are motivated by the work of Frohman, Gelca, Przytycki, Sikora
and others on the Kauffman bracket skein module, and its relation to the 
colored Jones function. Let us recall in brief these beautiful ideas.

For a manifold $N$, of any dimension,
possibly with nonempty boundary, let $\calS_q(N)$ denote 
the {\em Kauffman bracket skein 
module}, which is an $\R$-module (where $\R$ is defined in \eqref{eq.R})
generated by the isotopy classes of 
{\em framed unoriented} links in $N$ (including the empty one), modulo the 
relations of Figure \ref{f.relations}.

\begin{figure}[htpb]
$$ 
\psdraw{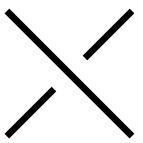}{0.5in}-q^{1/2} \psdraw{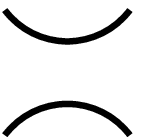}{0.5in} -q^{-1/2} \psdraw{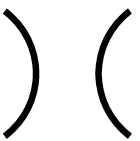}{0.5in},
\hspace{1cm}
\psdraw{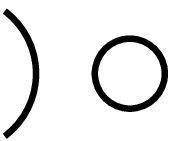}{0.5in} =(-q^{}-q^{-1}) \psdraw{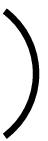}{0.15in}
$$
\caption{The relations of the Kauffman bracket skein module}\lbl{f.relations}
\end{figure}

\begin{remark} 
\lbl{rem.differs}
Our notation differs slightly from Gelca's et al \cite{FGL,Ge1,Ge2}. 
Gelca's $t^2$ equals to $q$ and further, Gelca's $(-1)^n J_n$ is our $J_{n+1}$
used below.
\end{remark}

Let us recall some elementary facts of skein theory, reminiscent of TQFT:
\begin{itemize}
\item[\bf{Fact 1:}]
If $N=N' \times I$, where $N'$ is a closed manifold, then
$\calS_q(N)$ is an algebra. 
\item[\bf{Fact 2:}]
If $N$ is a manifold with boundary $\partial N$,
then $\calS_q(N)$ is a module over the algebra $\calS_q(\partial N \times I)$.
\item[\bf{Fact 3:}]
If $N=N_1 \cup_Y N_2$ is the union of $N_1$ and $N_2$ along their common
boundary $Y$, then there is a map:
$$ \la \, , \, \ra: 
\calS_q(N_1) \otimes_{\calS_q(\partial Y \times I)}
\calS_q(N_1) \longto \calS_q(N)
$$
\end{itemize} 

We will apply the previous discussion in the following situation. Let $K$ 
denote a knot in a homology sphere $N$, and let $M$ denote the complement
of a thickening of $K$. Then, using the abbreviation 
$\calS_q(\BT):=\calS_q(\BT^2 \times I)$, Gelca and Frohman introduced in 
\cite{FG,FGL} the quantum peripheral and orthogonal ideals of $K$ (the latter 
was called {\em formal} in \cite[Sec.5]{FGL}):
\begin{definition}
\lbl{def.PO}
\rm{(a)} We define the {\em quantum peripheral ideal} 
$\calP(K)$ of $K$ to be the {\em annihilator} of the action of 
$\calS_q(\BT)$ on $\calS_q(M)$.
In other words,
$$
\calP(K)=\{ P \in \calS_q(\BT)| P. \emptyset =0 \}
$$
\rm{(b)} We define the {\em quantum orthogonal ideal} of a knot $K$ to be
$$
\calO(K)=\{v\in \calS_q(\BT)| \,  
\la \calS_q(S^1\times D^2)v, \emptyset \ra=0\}.
$$
\end{definition}

Note that $\calO(K)$ and $\calP(K)$ are left ideals in $\calS_q(\BT)$ 
and that $\calP(K) \subset \calO(K)$. Unfortunately, the quantum peripheral
and orthogonal ideals of a knot do not seem to be algorithmically computable
objects.

To understand the quantum peripheral and orthogonal ideals requires a better
description of the ring $\calS_q(\BT)$, 
and its module $\calS_q(S^1\times D^2)$.

The skein module $\calS_q(F \times I)$ is well-studied for a closed surface
$F$. It is a free $\R$-module on the set of free homotopy classes
of finite (possibly empty) collections of disjoint unoriented curves in $F$ 
without contractible components.
In particular, for a 2-torus $\BT$, $\calS_q(\BT)$ is the quotient of
the free $\R$-module on the set $\{ (a,b) \, | a, b \in \BZ\}$ 
modulo the relations $(a,b)=(-a,-b)$.
The multiplicative structure of $\calS_q(\BT)$
is well-known, and related to the {\em even trigonometric polynomials} of the
{\em quantum torus}, \cite{FG}. Let us recall this description due
to Gelca and Frohman. Consider the ring involution 
given by
\begin{equation}
\lbl{eq.tau}
\tau: \C \longto \C \hspace{1cm} E^a Q^b \longto E^{-a} Q^{-b}
\end{equation}
and let $\C^{\BZ_2}$ denote the invariant subring of $\C$.
Frohman-Gelca \cite{FG} prove that 
\begin{itemize}
\item[\bf{Fact 4:}]
The map
\begin{equation}
\lbl{eq.Phi}
\Phi: \calS_q(\BT) \longrightarrow \C^{\BZ_2}
\end{equation}
given by
$$
(a,b) \longto (-1)^{a+b} q^{-ab/2}(E^aQ^b+E^{-a}Q^{-b}).
$$
when $(a,b) \neq (0,0)$ and $\Phi(0,0)=1$, 
is an isomorphism of rings.
\end{itemize}

Thus, to a knot one can associate three ideals: the recursion ideal in
$\C$ and the quantum peripheral and the orthogonal ideal in $\calS_q(\BT)$. 
The next
theorem explains the relation between the recursion and quantum orthogonal 
ideals.

\begin{theorem}
\lbl{thm.2}
\rm{(a)} We have:
$$\Phi(\calO)=\C^{\BZ_2} \cap \calI.
$$
In particular, the colored Jones function of a knot determines its
quantum orthogonal ideal.
\newline
\rm{(b)} $\calI$ is invariant under the ring involution $\tau$.
\end{theorem}

The proof of the above theorem proves the following corollary which
compares the {\em orthogonality relations} of Gelca \cite[Sec.3]{Ge1} 
with the recursion relations given here:

\begin{corollary}
\lbl{cor.thm2}
Fix an element $x$ of the quantum orthogonal ideal of a knot. 
The orthogonality 
relation for the colored Jones function is $(E-E^{-1})x J=0$. On the other
hand, Theorem \ref{thm.2} implies that $xJ=0$. It follows that a $(d+2)$-term
recursion relation for the colored Jones function given by Gelca is
implied by a $d$-term recursion relation.
\end{corollary}

\begin{remark}
\lbl{rem.error}
Frohman and Gelca claimed that the quantum peripheral ideal is nonzero, 
\cite[Prop.8]{FGL} and \cite[Proof of Cor.1]{Ge1}, by specializing at $q=1$
and using the fact that the classical peripheral ideal is nonzero.
At the time of the \cite{FGL} paper, the nontriviality of the classical 
peripheral ideal was known for hyperbolic knots. Later on, the nontriviality
was shown by N. Dunfield and the author for all knots in $S^3$, \cite{DG}.

Combined with
our Theorem \ref{thm.2} below, the surjection of the specialization
map (that sets $q=1$) would prove Theorem \ref{thm.GL}.
Unfortunately, there is an error in the argument of Frohman and Gelca.
Instead, the surjection of the specialization map 
became known as the AJ Conjecture, formulated by the author in \cite{Ga2}. 
For a state-of-the art knowledge on the AJ Conjecture, see \cite{Hi}
and especially \cite{Le}. The latter gives a friendly discussion of the
classical and quantum peripheral and orthogonal ideals of a knot.
\end{remark}

\subsection{Converting difference into differential equations}
\lbl{sec.loop}



We now have all the ingredients to translates difference equations for 
$\{J_{K,n}\}$ to differential equations for $\{Q_{K,k}\}$.

\begin{theorem}
\lbl{thm.3}
\rm{(a)} Theorem \ref{thm.GL} implies a hierarchy of ODEs for
$\{Q_{K,k}\}$. More precisely, for every knot there exist a lower diagonal
matrix of infinite size
\begin{equation}
\lbl{eq.hierarchy}
D=\begin{pmatrix}
D_0 & 0 & 0 & \dots \\
D_1 & D_0 & 0 & \ddots \\
D_2 & D_1 & D_0 & \ddots \\
\dots  & \dots & \dots & \ddots 
\end{pmatrix}
\end{equation}
such that $D_i \in A_1$, $0 \neq D_0$ 
and $D Q_K=0$ where $Q_K=(Q_{K,1}, Q_{K,2}, \dots)^T$.
\newline
\rm{(b)} The above hierarchy uniquely determines the sequence $\{Q_{K,k}\}$
up to a finite number of initial conditions $\{ \frac{d^j}{du^j}|_{u=j}
Q_{K,k}(u)\}$ 
for $0 \leq j \leq \deg(D_0)$. 
\end{theorem}

Notice that the hierarchy \eqref{eq.hierarchy} depends on 
a linear $q$-difference equation for $\{J_{K,n}\}$. 
The degree of $D_0$ can be computed
from a linear $q$-difference equation for $\{J_{K,n}\}$, see \eqref{eq.degd}
below.

\subsection{Regular knots}
\lbl{sub.regular}

In this section we introduce the notion of a regular knot, and explain
the importance of this class of knots. The smallest degree of a nontrivial
differential operator is $1$.

\begin{definition}
\lbl{def.regular}
A knot $K$ is {\em regular} if $J_K$ satisfies a $q$-difference equation so
that $\deg(D_0)=1$.
\end{definition}

The explicit formulas of \cite{GL1} imply that the knots $3_1$ and $4_1$
are regular. X. Sun and the author will prove in a forthcoming paper that
all {\em twist knots} are regular. 

Among other reasons, regularity is important because of the following.

\begin{corollary}
\lbl{cor.regular}
If a knot $K$ is regular, then $\Jrat_K(q,u)$ (and thus, also $J_K$) is
uniquely determined by the hierarchy \eqref{eq.hierarchy} and the
initial condition
\begin{equation}
\lbl{eq.kashaev0}
\Jrat_K(q,1)=\sum_{k=0}^\infty Q_{K,k}(1)(q-1)^k \in \BQ[[q-1]]
\end{equation}
\end{corollary}

The power series invariant $\Jrat_K(q,1)$ is a disguised form of the
{\em Kashaev invariant} of a knot, which plays a prominent role in the
Volume Conjecture, due to Kashaev, and H. and J. Murakami; see \cite{Kv}
and \cite{MM}. Recall that the {\em Volume Conjecture} states that for a 
hyperbolic knot $K$ we have
$$
\lim_{n \to \infty} \frac{1}{n} \log|J_{K,n}(e^{2 \pi i/n})|
=\frac{1}{2 \pi} \text{vol}(K)
$$
where $\text{vol}(K)$ is the volume of $K$; \cite{Th}. In \cite{HL}, 
TTQ. Le and H. Vu reformulated the sequence 
$\{J_{K,n}(e^{2 \pi i/n})\}$ in terms of 
evaluation of a single function
$$
\kappa_K(q) \in \Lhat:=\invlim_j \BZ[q^{\pm}]/((1-q)(1-q^2) \dots (1-q^j)).
$$
More precisely, Le constructed an element $\kappa_K(q)$ in the above ring
so that for all $n \geq 1$ we have:
$$
\kappa_K(e^{2 \pi i/n})=J_{K,n}(e^{2 \pi i/n}).
$$
There is a Taylor series map:
$$
\text{Taylor}: \Lhat \longto \BZ[[q-1]].
$$
Then, we have:
\begin{equation}
\lbl{eq.important}
\Jrat_K(q,1)=\text{Taylor}(\kappa_K(q)).
\end{equation}
The reader may deduce a proof of the above equation from \cite[Sec.3]{GL2}.
Thus, for regular knots, the Kashaev invariant together with the ODE
hierarchy \eqref{eq.hierarchy} {\em uniquely determines} the colored Jones
function $\{J_K\}$, and its growth-rate in the Generalized Volume
Conjecture.

\subsection{Questions}
\lbl{sub.questions}

The above hierarchy is reminiscent of {\em matrix models} discussed
in physics. See for example \cite{DV} and Question \ref{que.1} below.

Let us mention that the conversion of difference into differential equations
can actually be interpreted as an application of the {\em WKB method} for the
linear $q$-difference equation. This remark, and its implications to
asymptotics of the colored Jones function is explained in a later publication,
\cite{GG}.

Let us end with some questions.

\begin{question}
\lbl{que.1}
Is there a {\em physical meaning} to the recursion relations of the
colored Jones function, and in particular of the hierarchy of ODE which
is satisfied by its loop expansion? Differential equations often hint
at a hidden {\em matrix model}, or an {\em M-theory} explanation.
\end{question}

\begin{question}
\lbl{que.2}
The hierarchy of ODEs that appear in Theorem \ref{thm.3} also appears,
under the name of {\em semi-pfaffian chain}, in complexity
questions of real algebraic geometry. We thank S. Basu for pointing
this out to us. For a reference, see \cite{GV}. Is this a coincidence?
\end{question}

\begin{question}
\lbl{que.3}
In \cite{GK} Kricker and the author constructed a rational form
$\Zrat$
of the Kontsevich integral of a knot. As was explained in \cite{Ga},
this rational form becomes the loop expansion of the colored Jones
function, on the level of Lie algebras. In \cite{GL1} it is shown that
the $\fg$-colored Jones function for any simple Lie algebra
$\fg$ is $q$-holonomic. Holonomicity gives rise to a ring $\C_{q,\fg}$
with an action of the Weyl group $W$ of $\fg$. The ring $\C_{q,\fg}$
specializes to the coordinate ring of the $\fg$-character variety of the 
torus, introduced by Przytycki-Sikora, at least in case $\fg=\mathfrak{sl}_n$;
see for example \cite{PS} and \cite{Si}.

Is there a Kauffman bracket
skein theory $\calS_{q,\fg}$ that depends on $\fg$, in such a way that we have 
an ring isomorphism:
\begin{equation*}
\Phi_{\fg}: \calS_{q,\fg}(\BT) \longrightarrow \C^{W}_{\fg}?
\end{equation*}
If so, one could define the $\fg$-quantum peripheral and orthogonal ideals
of a knot, and ask whether the analogue of Theorem \ref{thm.2} holds:
$$
\Phi_{\fg}(\calO_{\fg})=\C_{\fg}^W \cap \calI_{\fg}?
$$
In addition, one may ask for an analogue of Theorem \ref{thm.3} for simple
Lie algebras $\fg$. Notice, however, that this analogue is going to be 
a hierarchy of PDEs for functions of as many variables as the rank of
the Lie algebra.
\end{question}

\begin{question}
\lbl{que.5}
Theorem \ref{thm.2} proves that $\calI$ is an ideal in $\C^{\BZ_2}$
invariant under the ring involution $\tau$. Is it true that
$\calI$ is generated by its $\BZ_2$-invariant part? In other words, is it true
that $\calI=\C (\C^{\BZ_2} \cap \calI)$? 
\end{question}

\subsection{Acknowledgement}
An early version of this paper was announced in the JAMI 2003 meeting in
Hohns Hopkins. We wish to thank Jack Morava for the invitation.
We wish to thank TTQ. Le and D. Zeilberger for stimulating 
conversations, and especially A. Sikora for enlightening conversations
on skein theory and for explaining to us 
beautiful work of Frohman and Gelca on the Kauffman 
bracket skein module.

\section{Proof of Theorem \ref{thm.2}}
\lbl{sec.related}

Let us begin by discussing the recursion relation for the colored Jones
function which is obtained by a nonzero element of the quantum orthogonal ideal
of a knot. This uses work of Gelca, which
we will quote here. For proofs, we refer the reader to \cite{Ge1}.

The problem is to understand the right action of the
ring $\calS_q(\BT)$ on the skein module $\calS_q(S^1 \times D^2)$.

To begin with, the skein module $\calS_q(S^1 \times D^2)$ can be identified
with the polynomial ring $\R[\a]$, where $\a$ is a longitudonal curve
in the solid torus, and $\R=\BZ[q^{\pm 1/2}]$.
Rather than using the $\R$-basis for  $\calS_q(S^1 \times D^2)$ given by
$\{\a^n\}_n$, Gelca uses the basis given by $\{T_n(\a)\}_n$, where
$\{T_n\}$ is a sequence of Chebytchev-like polynomials defined by 
$T_0(x)=2$, $T_1(x)=x$ and $T_{n+1}(x)=xT_n(x)-T_{n-1}(x)$.

Recall that the ring $\calS_q(\BT)$ is generated by symbols $(a,b)$ for 
integers $a$ and $b$ and relations $(a,b)=(-a,-b)$. 
Gelca \cite[Lemma 1]{Ge1} describes the right action of $\calS_q(\BT)$ on 
$\calS_q(S^1 \times D^2)$ as follows: 
\begin{eqnarray}
\lbl{eq.action}
T_n(\a) \cdot (a,b) &=&
q^{ab/2}(-1)^b
\left( q^{nb} \left[
      q^{b} S_{n+a}(\a) 
      -q^{-b} S_{n+a-2}(\a) \right] \right. \\
\notag
& & + \left. q^{-nb} \left[
      -q^{b} S_{n-a-2}(\a) 
      +q^{-b} S_{n-a}(\a) \right] \right) ,
\end{eqnarray}
where $\{S_n\}$ is a sequence of Chebytchev polynomials defined by 
$S_0(x)=1$, $S_1(x)=x$ and $S_{n+1}(x)=xS_n(x)-S_{n-1}(x)$.

Consider an element $x=\sum_{a,b} c_{a,b} (a,b)$ of the quantum orthogonal 
ideal of 
a knot and recall the pairing $\la \, , \, \ra$ from Fact 3.
Since the (shifted) $n$th colored Jones polynomial of a knot is given by
$J_n=(-1)^{n-1} \la S_{n-1}(\a), \emptyset \ra$, Equation \eqref{eq.action} 
implies the recursion relation 
\begin{eqnarray}
\lbl{eq.rec}
0 &=& \sum_{a,b} c_{a,b} q^{ab/2} (-1)^{a+b} \left( q^{nb} \left[
      q^{b} J_{n+a+1}(K) 
      -q^{-b} J_{n+a-1}(K) \right] \right. \\
\notag
& & + \left. q^{-nb} \left[
      -q^{b} J_{n-a-1}(K) 
      +q^{-b} J_{n-a+1}(K) \right] \right)
\end{eqnarray}
for the colored Jones function
corresponding to an element $\sum_{a,b} c_{a,b}(a,b)$ in the quantum orthogonal
ideal.

Let us write the above recursion relation in operator form.
Recall that the operators $E$ and $Q$ act on the discrete function $J$ by
$(EJ)(n)=J(n+1)$ and $QJ(n)=q^n J(n)$, and satisfy the commutation
relation $EQ=qQE$. Then, Equation \eqref{eq.rec} becomes:
$$
0= \sum_{a,b} c_{a,b} q^{ab/2} (-1)^{a+b} \left( q^b Q^{b} E^{a+1}
-q^{-b} Q^{b} E^{a-1} - q^{b} Q^{-b} E^{-a-1} 
+ q^{-b} Q^{-b} E^{-a+1} \right)J.
$$
Using the commutation relation $E^k Q^l = q^{kl} Q^l E^k$ for integers
$k,l$ and moving the $E$'s on the left and the $Q$'s on the right,
we obtain that
$$
0 = (E-E^{-1}) \sum_{a,b} c_{a,b} q^{-ab/2} (E^a Q^b + E^{-a} Q^{-b})J.
$$
Recall the isomorphism $\Phi$ of Equation \eqref{eq.Phi}. Using this
isomorphism, our discussion so far implies that $x \in \C^{\BZ_2}$ is
an element of the quantum orthogonal ideal of a knot iff $(E-E^{-1})x $ lies
in the recursion ideal. It remains to show that for every $x \in \C^{\BZ_2}$,
$(E-E^{-1})x \in \calI$
iff $x \in \calI$. One direction is obvious since $\calI$ is a left
ideal. For the opposite direction, consider $x \in \C^{\BZ_2}$, and
let $y=(E-E^{-1})x$ and $f=xJ$. Assume that $y \in \calI$. We need to show
that $x \in \calI$; in other words that $f=0$.

We have $(E^2-I)f=E(E-E^{-1})f=0$, which implies that
\begin{equation}
\lbl{eq.f1}
f(n+2)=f(n)
\end{equation}
for all $n \in \BZ$.

Recall the symmetry relation $J_n+J_{-n}=0$ for the colored Jones function. 
In order to write it in operator form, consider the operator
$S$ that acts on a discrete function $f$ by $(Sf)(n)=f(-n)$.
Then, $(S+I)J=0$.

It is easy to see that 
\begin{equation}
\lbl{eq.S}
SE=E^{-1}S \hspace{1cm} SQ=Q^{-1}S.
\end{equation}
Since $\C^{\BZ_2}$
is generated by $E^aQ^b+E^{-a}Q^{-b}$, it follows that
$S$ commutes with every element of $\C^{\BZ_2}$; in particular
$Sx=xS$, and thus $(S+I)f=(S+I)xJ=x(S+I)J=0$. In other words, 
\begin{equation}
\lbl{eq.f2}
f(n)+f(-n)=0
\end{equation}
for all $n$. Equations \eqref{eq.f1} and \eqref{eq.f2} imply that
$f(2n)=f(0)$, $f(2n+1)=f(1)$, $f(0)=f(1)=0$. Thus, $f=0$. This completes
part (a) of Theorem \ref{thm.2}.

For part (b), consider $x \in \calI$ and recall the involution $\tau$
of \eqref{eq.tau}. Then, we have $x=x_+ + x_-$
where $x_{\pm}=1/2(x \pm \tau(x)) \in \C_\pm$, where $\C_\pm$ is generated by
$E^a Q^b \pm E^{-a} Q^{-b}$. \eqref{eq.S} implies that $Sx_+=x_+S$
and $S x_-=x_-S$. 
Now, we have
$$
0=xJ=x(-J)=xSJ=(x_+ +x_-)SJ=S(x_+ -x_-)J.
$$
Since $S^2=I$, it follows that $(x_+ -x_-)J=0$. This, together with
$0=(x_+ +x_-)J$, implies that $x_\pm J=0$. In other words, $x_\pm \in \calI$.
Since $\tau(x)=x_+-x_-$, it follows that $\calI$ is invariant under $\tau$.
\qed

The proof of Theorem \ref{thm.2} also proves Corollary \ref{cor.thm2}.

\begin{example}
\lbl{ex.1}
In \cite{Ge2} Gelca computes that the following element
$$
(1,-2k-3)-q^{-4}(1,-2k+1)+q^{(2k-5)/2}(0,2k+3)-q^{(2k-1)/2}(0,2k-1)
$$
lies in the quantum peripheral (and thus quantum orthogonal) ideal of the 
left handed 
$(2,2k+1)$ torus knot. Using the isomorphism $\Phi$ and Theorem \ref{thm.2} 
and a simple calculation, it follows that the following element
$$
-q^2(EQ^{-2k-3}+E^{-1}Q^{2k+3}) +q^{-4}(EQ^{-2k+1}+E^{-1}Q^{2k-1})
+q^{-2}(Q^{2k+3}+Q^{-2k-3}) -(Q^{2k-1}+Q^{-2k+1})
$$
lies in the recursion ideal of the left handed $(2,2k+1)$ torus knot. 
This element
gives rise to a $3$-term recursion relation for the colored Jones function.
\end{example}

\section{Proof of Theorem \ref{thm.3}}
\lbl{sec.proof3}

\subsection{Three $q$-difference rings}
\lbl{sub.aux}

In this section we consider some auxiliary rings and their associated
$\C$-module structure.

The colored Jones function is a sequence of Laurent polynomials, in other
words an element of the ring $\BZ[q^{\pm/2}]^{\BZ}$. The ring 
$\BZ[q^{\pm/2}]^{\BZ}$ is a $\C$-module via the action \eqref{eq.EQaction1}.
In other words, for $f \in \BZ[q^{\pm/2}]^{\BZ}$, we define:
$$
(Q f)_n(q) = q^n f_n(q), \quad
(E f)_n(q) = f_{n+1}(q).
$$
Consider the ring $\BQ'(u)[[q-1]]$ from Section \ref{sub.goal}. 
It is  also a $\C$-module, where for $f(q,u) \in \BQ'(u)[[q-1]]$, we define:
$$
(Q f)(q,u) = u f(q,u), \quad
(E f)(q,u) = f(q,qu).
$$
Consider in addition the ring $\BQ[[q-1]]^{\BZ}$. 
It is  a $\C$-module, where for $(f_n(q)) \in \BQ[[q-1]]$, we define:
$$
(Q f)_n(q) = q^n f_n(q), \quad
(E f)_n(q) = f_{n+1}(q).
$$
In the language of {\em differential Galois theory}, $\BZ[q^{\pm/2}]^{\BZ}$,
$\BQ'(u)[[q-1]]$ and $\BQ[[q-1]]^{\BZ}$ are {\em $q$-difference rings},
see for example \cite{vPS}.
There is a ring homomorphism:
$$
\Ev: \BQ'(u)[[q-1]] \longto \BQ[[q-1]]^{\BZ}
$$
given by 
$$
f(q,u) \mapsto (f(q,q^n)).
$$
which respects the $\C$-module structure. In other words, for $f(q,u) \in 
\BQ'(u)[[q-1]]$, we have:
$$
\Ev(Qf)=Q\Ev(f), \quad \Ev(Ef)=E\Ev(f).
$$
It is not hard to see that $\Ev$ map is injective. 
Equation \eqref{eq.Jrat} then states that for every knot $K$ we have:
$$ 
\Ev(\Jrat_K)=J_K/J_{\text{unknot}}.
$$
Notice further that since $J_K$ is holonomic and 
$J_{\text{unknot},n}(q)=[n]$ is {\em closed-form} (that is,
$[n+1] \in \BQ(q^{n/2}, q^{1/2}$, it follows that 
$J_K/J_{\text{unknot}}$ is holonomic, too. In the proof of Theorem
\ref{thm.3} below, $J$ and $\Jrat$ will stand for $J_K/J_{\text{unknot}}$ 
and $\Jrat_K$ respectively.

\subsection{Proof of Theorem \ref{thm.3}}
\lbl{sub.proof3}

Consider a function 
\begin{equation}
\lbl{eq.myrat}
\Jrat(q,u) =\sum_{k=0}^\infty Q_k(u) (q-1)^k \in \BQ'(u)[[q-1]]
\end{equation}
such that
$J:=\Ev(\Jrat) \in \BQ[[q-1]]^{\BZ}$ is holonomic. Thus, 
$X J=0$ where $X =\sum_{a,b} c_{a,b}(q) E^a Q^b \in \C$, and 
the sum is finite. In other words,
we have:
$$
0 = \sum_{a,b} c_{a,b}(q) q^{(n+a)b} J_{n+a}(q).
$$
Since $\Psi$ is a map of $\C$-modules, it follows that
$$
0 = \Ev\left( \sum_{a,b} \ti c_{a,b}(q) Q^b E^a \Jrat\right),
$$
where
$$
\ti c_{a,b}(q)=c_{a,b}(q) q^{ab}.
$$
Since $\Psi$ is an injection, it follows that
$$
0 =  \sum_{a,b} \ti c_{a,b}(q) Q^b E^a \Jrat(q,u).
$$
Using Equation \eqref{eq.myrat} and interchanging the order of summation,
we obtain that
\begin{eqnarray}
\notag 
0 &=& \sum_{a,b} \ti c_{a,b}(q) u^b  \Jrat(q,q^a u) \\
\lbl{eq.XQ} 
&=& \sum_{k=0}^\infty (q-1)^k \calX Q_k,
\end{eqnarray}
where 
$$
\calX: \BQ'(u) \longto \BQ'(u)[[q-1]]
$$ 
is the operator defined by
\begin{eqnarray*}
f & \mapsto & \calX f=\sum_{a,b} c_{a,b}(q) u^b q^{ab} f(uq^a) \\
& = & \sum_{a,b} \ti c_{a,b}(q) s^b f(uq^a)
\end{eqnarray*}
Let us define
$$
P(\l,u,q)=\sum_{a,b} \ti c_{a,b}(q) u^b \l^a \in 
\BZ[\l^{\pm},u^{\pm},q^{\pm}].
$$
In the language of $q$-difference equations, $P(\l,u,q)$ is often called
the {\em characteristic polynomial} of $X$.
Let us denote by $\la f \ra_m$ the coefficient of $(q-1)^m$ in a power 
series $f$. 
Applying $\la \cdot \ra_m$ to \eqref{eq.XQ}, it follows that for 
all $m \geq 0$ we have
\begin{equation}
\lbl{eq.DE2}
0 = \la \calX Q_0 \ra_m +
\la \calX Q_1 \ra_{m-1} 
 \dots + \la \calX Q_m \ra_0
\end{equation}
The chain rule implies that 
\begin{equation}
\lbl{eq.D'}
\la \calX f \ra_m=\calD_m f
\end{equation}
for some differential
operator $\calD_m$ with polynomial coefficients
of degree at most $m$. For example, we have:
\begin{eqnarray*}
(\calD_0 f)(u) &=& P(1,u,1) f(u) \\
(\calD_1 f)(u) &=& P_q(1,u,1) f(u) + P_{\l}(1,u,1) u f'(u) \\
(\calD_2 f)(u) &=& P_{qq}(1,u,1) f(u) + P_{\l q}(1,u,1)u f'(u) + 
P_{\l\l}(1,u,1)u^2 f''(u) + (P_{\l\l}(1,u,1)-P_{\l}(1,u,1)) u f'(u), 
\end{eqnarray*}
where the subscripts $._q$ and $._{\l}$ denote $\pt/\pt q$ and $\l \pt/\pt \l$
respectively, and the superscript denotes derivative with respect to $u$.

Thus, we obtain that for all $m \geq 0$
\begin{equation}
\lbl{eq.DE}
0 = \calD_m Q_0 + \calD_{m-1}Q_1  \dots + \calD_0 Q_m. 
\end{equation}
We will show shortly that $\calD_m \neq 0$ for some $m$. Assuming this, let
$l=\min\{ m | \calD_m \neq 0\}$ and define $D_m=\calD_{l+m}$.
Equation \eqref{eq.DE} for $l+m$ implies that
$$
0= D_{m+l} Q_0  + D_{m+l-1} Q_1 + \dots D_0 Q_{m+l}.
$$
In other words,
$$
\begin{pmatrix}
D_0 & 0 & 0 & \dots \\
D_1 & D_0 & 0 & \ddots \\
D_2 & D_1 & D_0 & \ddots \\
\dots  & \dots & \dots & \ddots 
\end{pmatrix} 
\begin{pmatrix}
Q_0 \\ Q_1 \\ Q_2 \\ Q_3 \\ \vdots 
\end{pmatrix} =
\begin{pmatrix}
0 \\ 0 \\ 0 \\ 0 \\ \vdots 
\end{pmatrix}
$$
as needed. 

It remains to prove that $\calD_m \neq 0$ for some $m$. 
The definition of $\calD_m$ 
implies easily that $\calD_k=0$ for $k \leq m$ if and only if for all 
multiindices
$I=(i_1,\dots, i_k)$ with $i_j \in \{q, \l \}$ and $k \leq m$, we have:
\begin{equation}
\lbl{eq.PI}
P_I(1,u,1)=0.
\end{equation}
Since $P(\l,u,q)$ is a Laurent polynomial, 
if $\calD_m=0$ for all $m$, then $P(\l,u,q)=0$. 

{\em Additively},
there is an $\BZ$-linear isomorphism $\C \leftrightarrow 
\BZ[\l^{\pm},u^{\pm 1},q^{\pm 1}]$, given by $Q^b E^a \mapsto u^b v^a$,
and sends $X \leftrightarrow P(\l,u,q)$. Thus, $X=0$ a contradiction
to our hypothesis.
This concludes part (a) of Theorem \ref{thm.3}. Part (b) follows
from Lemma \ref{lem.ode} below.
\qed

Let us finally compute $l$ and the degree $d$ of $D_0\calD_l$ from
$X$. Notice that $d$ equals to the number of initial
conditions needed to determine the sequence $\{Q_k\}$ from the ODE hierarchy
\eqref{eq.hierarchy}.

For natural numbers $n,m$ with $n \geq m$, let us denote by
$ I(n,m)=(\l,\dots,\l,q,\dots,q)$ the multiindex of length $|I(n,m)|=n$
where $\l$ appears $m$ times and $q$ appears $n-m$ times.  

Equation \eqref{eq.PI} implies that
\begin{equation}
\lbl{eq.degl}
l=\min \{n \, | P_{I(n,m)}(1,u,1) \neq 0 \quad \text{for some} \quad m \}
\end{equation}
and
\begin{equation}
\lbl{eq.degd}
d=\min \{m \, | P_{I(l,m)}(1,u,1) \neq 0  \}.
\end{equation}

\begin{lemma}
\lbl{lem.ode}
If $a_i,c \in \BC[u^{\pm 1}]$, $a_n \neq 0$, the ODE 
$$ 
a_n f^{(n)} + a_{n-1} f^{(n-1)} + \dots a_0 f=c
$$
has at most one solution which is a rational function with fixed 
initial condition for $f^{(k)}(x_0)$ for $k=0, \dots, n-1$, where
$a_n(x_0) \neq 0$.
\end{lemma}

\begin{proof}
Consider the set of real numbers $u$ such that $a_0(u) \neq 0$.
It is a finite set of open intervals. Uniqueness of the solution (modulo
initial conditions) is well-known. Since $f$ is a rational function,
it is uniquely determined by its restriction on an open interval. 
The result follows.
\end{proof}

\ifx\undefined\bysame
        \newcommand{\bysame}{\leavevmode\hbox
to3em{\hrulefill}\,}
\fi


\begin{thebibliography}{[EMSS]}

\bibitem[B1]{B1} I.N. Bernstein,
        {\em Modules over a ring of differential operators. An investigation 
        of the fundamental solutions of equations with constant coefficients},
        Functional Analysis and its applications {\bf 5} 1971 1--16, 
        English translation: 89--101.

\bibitem[B2]{B2} \bysame,
        {\em The analytic continuation of generalized functions with respect
        to a parameter},
        Functional Analysis and its applications {\bf 6} 1972 26--40,
        English translation: 273--285.

\bibitem[Bj]{Bj} J-E. Bj\"ork,
        {\em Rings of differential operators},
        North-Holland, 1971.

\bibitem[DV]{DV} R. Dijkgraaf and C. Vafa,
        {\em The geometry of matrix models},
        preprint 2002 {\tt hep-th/0207106}.

\bibitem[DG]{DG} N. Dunfield and S. Garoufalidis,
        {\em Non-triviality of the $A$-polynomial for knots in $S^3$},
        Algebr. Geom. Topol. {\bf 4} (2004) 1145-1153.

\bibitem[FG]{FG} C. Frohman and R. Gelca, 
        {\em Skein modules and the noncommutative torus},
        Trans. Amer. Math. Soc. {\bf 352} (2000) 4877--4888.

\bibitem[FGL]{FGL} \bysame, \bysame and W. Lofaro,
        {\em The A-polynomial from the noncommutative viewpoint},
        Trans. Amer. Math. Soc. {\bf 354} (2002)  735--747.

\bibitem[GV]{GV} A. Gabrielov and N. Vorobjov,
        {\em Complexity of stratifications of semi-Pfaffian sets},
        Discrete Comput. Geom. {\bf 14} (1995)  71--91.

\bibitem[GK]{GK} S. Garoufalidis and A. Kricker,
        {\em A rational noncommutative invariant of boundary links}, 
        Geom. and Topology {\bf 8} (2004) 115--204.

\bibitem[Ga]{Ga} \bysame,
        {\em Rationality: From Lie algebras to Lie groups},
        preprint 2002 {\tt math.GT/0201056}.

\bibitem[GL]{GL1} \bysame and T.T.Q. Le,
        {\em The colored Jones function is $q$-holonomic}
        Geom. and Topology {\bf 9} (2005) 1253--1293.

\bibitem[GL]{GL2} \bysame and \bysame,
        {\em Asymptotics of the colored Jones function of a knot},
        preprint 2005.

\bibitem[Ga2]{Ga2} \bysame,
        {\em On the characteristic and deformation varieties of a knot},
        Proceedings of the CassonFest,
        Geometry and Topology Monographs {\bf 7} (2004) 291--309.

\bibitem[GG]{GG} \bysame and J. Geronimo,
        {\em  Asymptotics of $q$-difference equations},
        Contemporary Math. AMS {\bf 416} (2006) 83--114.

\bibitem[Ge1]{Ge1} R. Gelca,
        {\em On the relation between the $A$-polynomial and the Jones 
        polynomial}, 
        Proc. Amer. Math. Soc. {\bf 130} (2002)  1235--1241.

\bibitem[Ge2]{Ge2} \bysame,
        {\em The noncommutative A-ideal of a (2,2p+1)-torus knot determines 
        its Jones polynomial},
        J. Knot Theory Ramifications  {\bf 12}  (2003) 187--201. 

\bibitem[Hi]{Hi} K. Hikami,
        {\em Difference equation of the colored Jones polynomial for torus 
        knot},
        Internat. J. Math.  {\bf 15}  (2004) 959--965.

\bibitem[HL]{HL} V. Huynh and T.T.Q. Le,
        {\em On the Colored Jones Polynomial and the Kashaev invariant},
         Fundam. Prikl. Mat.  {\bf 11}  (2005) 57--78.

\bibitem[J]{J} V. Jones,
        {\em Hecke algebra representation of braid groups and link
        polynomials}, 
        Annals Math. {\bf 126} (1987)  335--388.

\bibitem[Kv]{Kv} R. Kashaev,
        {\em The hyperbolic volume of knots from the quantum dilogarithm},
        Modern Phys. Lett. A {\bf 39} (1997) 269--275.

\bibitem[Ka]{Ka} C. Kassel,
        {\em Quantum groups}, 
        GTM {\bf 155} Springer-Verlag (1995).

\bibitem[Le]{Le} T.T.Q. Le,
        {\em The Colored Jones Polynomial and the $A$-Polynomial of 
        Two-Bridge Knots},
        Advances Math. in press.

\bibitem[M]{M} Y. Manin,
        {\em Quantum groups and noncommutative geometry},
        Universit\'e de Montr\'eal, Centre de Recherches Math\'ematiques, 
        Montreal, QC, 1988.

\bibitem[MM]{MM} H. Murakami and J. Murakami,
        {\em The colored Jones polynomials and the simplicial volume of a 
        knot},
        Acta Math.  {\bf 186}  (2001) 85--104.

\bibitem[PS]{PS} J. Przytycki and A. Sikora, 
        {\em Skein algebra of a group},
        Knot theory,  297--306, Banach Center Publ., {\bf 42}, 
        Polish Acad. Sci., Warsaw, 1998.
 
\bibitem[PWZ]{PWZ} M. Petkov\v sek, H.S. Wilf and D.Zeilberger, 
        {\em $A=B$}, 
        A.K. Peters, Ltd., Wellesley, MA 1996.

\bibitem[Ro]{Ro} L. Rozansky,
        {\em The universal $R$-matrix, Burau representation, and the 
        Melvin-Morton expansion of the colored Jones polynomial},
        Adv. Math.  {\bf 134}  (1998) 1--31.
        A.K. Peters, Ltd., Wellesley, MA 1996.

\bibitem[Si]{Si} A. Sikora, 
        {\em ${\rm SL}\sb n$-character varieties as spaces of graphs},
        Trans. Amer. Math. Soc.  {\bf 353}  (2001) 2773--2804.

\bibitem[Th]{Th} W. Thurston, {\em The geometry and topology of 3-manifolds},
        1979 notes, available from MSRI.

\bibitem[Tu]{Tu} V. Turaev,
        {\em The Yang-Baxter equation and invariants of links},
        Inventiones Math. {\bf 92} (1988) 527--553.

\bibitem[vPS]{vPS} M. van der Put and M. Singer,
        {\em  Galois theory of difference equations},
        Lecture Notes in Mathematics, {\bf 1666} Springer-Verlag  1997.

\bibitem[WZ]{WZ} H. Wilf and D. Zeilberger,
        {\em An algorithmic proof theory for hypergeometric (ordinary and 
        $q$) multisum/integral identities}, 
        Inventiones Math. {\bf 108} (1992)  575--633.

\bibitem[Z]{Z} D. Zeilberger,
        {\em A holonomic systems approach to special functions identities},
        J. Comput. Appl. Math. {\bf 32} (1990) 321--368.

\end{thebibliography}
\end{document}